\documentclass[11pt]{amsart}
\usepackage[all]{xy}
\usepackage{amsmath}
\usepackage{amsfonts}
\usepackage{amssymb}
\usepackage{latexsym}
\usepackage{graphicx}
\usepackage{color}
\usepackage{amscd}
\usepackage{epsfig}
\usepackage{cite}

\newtheorem{problem}{Problem}[section]
\newtheorem{definition}[problem]{Definition}
\newtheorem{lemma}[problem]{Lemma}
\newtheorem{theorem}[problem]{Theorem}

\newtheorem{corollary}[problem]{Corollary}

\title{Twisted exponential sums of polynomials in one variable}
\author{Chunlei Liu}
\address{Department of Mathematics, Shanghai Jiao Tong
University, Shanghai 200240, P.R. China, E-mail: clliu@sjtu.edu.cn}
\author{Wenxin Liu}\address{School of Mathematical Sciences, Beijing Normal
University, Beijing 100875, P.R. China, E-mail:
wenxin8210@mail.bnu.edu.cn}
\thanks{This research is supported by NSFC Grant No.
10671015.\\
Dedicated to Yuan Wang on the occasion of his 80th
birthday}
\begin{document}
\maketitle

\begin{abstract}
The twisted $T$-adic exponential sum associated to a  polynomial in
one variable is studied.  An explicit arithmetic polygon in terms of
the highest two exponents of the polynomial is proved to be a lower
bound of the Newton polygon of the $C$-function of the twisted
T-adic exponential sum. This bound gives lower bounds for the Newton
polygon of the $L$-function of twisted $p$-power order exponential
sums.
\end{abstract}


\section {Introduction}

 Let $p$ be a prime number,  $q$ a power of
$p$, and $\mathbb{F}_q$ the finite field with $q$ elements. Let $W$
be the Witt ring scheme, $\mathbb{Z}_{q}=W(\mathbb{F}_q)$, and
$\mathbb{Q}_{q}=\mathbb{Z}_{q}[\frac{1}{p}]$.  Let $\mu_{q-1}$ be
the group of $(q-1)$-th roots of unity in $\mathbb{Z}_q$,
$\omega:x\mapsto\hat{x}$ the  Teichm\"{u}ller lifting from
$\mathbb{F}_q $  to $\mu_{q-1}$,  and $\chi=\omega^{-u}$ with
$u\in\mathbb{Z}^n/(q-1)$ a  character of $(\mathbb{F}_q^{\times})^n$
into $\mu_{q-1}$ .

Let $\triangle\supsetneq\{0\}$ be an integral convex polytope in
$\mathbb{R}^n$, and $I$ the set of vertices of $\triangle$ different
from the origin. Let
$$f(x)=\sum\limits_{u\in\triangle}(a_ux^u,0,0,\cdots)\in W(\mathbb{F}_q[x_1^{\pm1},x_2^{\pm1},\cdots,x_n^{\pm1}])\text{ with }
\prod_{u\in I}a_u\neq0,$$ where $x^u=x_1^{u_1}x_2^{u_2}\cdots
x_n^{u_n}$ if $u=(u_1,u_2,\cdots,u_n)\in \mathbb{Z}^n$.

\begin{definition}
For any positive integer $l$,  the sum
$$S_{f,\chi}(l,T)=\sum_{x\in(\mathbb{F}_{q^l}^{\times})^n}\chi({\rm Norm}_{\mathbb{F}_{q^l}/\mathbb{F}_q}(x))
(1+T)^{Tr_{\mathbb{Q}_{q^l}/\mathbb{Q}_p}(f(x))}\in\mathbb{Z}_q[[T]]$$
 is called a twisted T-adic exponential sum of
$f(x)$ . And the function

$$L_{f,\chi}(s,T)=\exp(\sum_{l=1}^{\infty}S_{f,\chi}(l,T)\frac{s^l}{l})$$
is called a  L-function of  twisted T-adic  exponential sums .
\end{definition}

We have
$$L_{f,\chi}(s,T)=\prod_{x\in|\mathbb{G}_m^n\otimes\mathbb{F}_q|}
\frac{1}{1-\chi({\rm
Norm}_{\mathbb{F}_{q^m}/\mathbb{F}_q}(x))(1+T)^{{\text
Fr}_x}s^{m}},$$ where $\mathbb{G}_m$ is the multiplicative group
$xy=1$,   $m=\deg(x)$ and ${\text
Fr}_x=Tr_{\mathbb{Q}_{q^m}/\mathbb{Q}_p}(f(x))$.

That Euler product formula gives
$$L_{f,\chi}(s,T)\in 1+s\mathbb{Z}_q[[T]][[s]].$$

The theory of $T$-adic exponential sums without twists was developed
by Liu-Wan \cite{LWn}, and the theory of twisted $T$-adic
exponential sums was developed by Liu \cite{Liu2}.

Let $m\geq1$, $\zeta_{p^m}$ a primitive $p^m$-th root of unity, and
$\pi_m=\zeta_{p^m}-1$. Then the specialization
 $L_{f,\chi}(s,\pi_m)$ is the $L$-function of twisted
$p$-power order exponential sums $S_{f,\chi}(l,\pi_m)$. These sums
were studied by Liu \cite{Liu}, with the $m=1$ case studied by
Adolphson-Sperber [AS,AS2], and if  $\chi$ is trivial, they were
studied by Liu-Wei \cite{LW}.

Define
$$C_{f,\chi}(s,T)=\exp(\sum_{l=1}^{\infty}-(q^l-1)^{-n}S_{f,\chi}(l,T)\frac{s^l}{l}).$$
Call it a $C$-function of twisted $T$-adic exponential sums.

 We
have
$$L_{f,\chi}(s,T)=\prod_{i=0}^n
C_{f,\chi}(q^is,T)^{(-1)^{n-i+1}{n\choose i}} ,$$ and
$$C_{f,\chi}(s,T)^{(-1)^{n-1}}=\prod_{j=0}^{+\infty}
L_{f,\chi}(q^js,T)^{{n+j-1\choose j}}.$$
 So we have
$$C_{f,\chi}(s,T)\in1+s\mathbb{Z}_q[[T]][[s]].$$

We  view $C_{f,\chi}(s, T)$ as a power series in the single variable
$s$ with coefficients in the $T$-adic complete field
$\mathbb{Q}_q((T))$. The $C$-function  $C_{f,\chi}(s,T)$ was shown
$T$-adic entire in $s$ by Liu \cite{Liu2}.
\\

Let $C(\triangle)$ be the cone generated by $\triangle$,
$M(\triangle)=C(\triangle)\cap\mathbb{Z}^n$, and $\deg_{\triangle}$
the degree function on $C(\triangle)$, which is $\mathbb{R}_{\geq
0}$ linear and takes the values $1$ on each co-dimension $1$ face
not containing $0$. Let $u\in\mathbb{Z}^n/(q-1)$, and
$$M_{u}(\triangle):=\frac{1}{q-1}(M(\triangle)\cap u).$$
\begin{definition}Let $b$ be the least positive integer such that $p^bu=u$. Order
elements of $\cup_{i=0}^{b-1}M_{p^iu}(\triangle)$ such that
$$\deg_{\triangle}(x_1)\leq\deg_{\triangle}(x_2)\leq\cdots.$$
A convex function on $\mathbb{R}_{\geq0}$ which is linear between
consecutive integers with initial value $0$ is called the infinite
u-twisted Hodge polygon of $\triangle$  if its slopes between
consecutive integers are the numbers
$$\frac{\deg_{\triangle}(x_{bi+1})+\deg_{\triangle}(x_{bi+2})+\cdots+\deg_{\triangle}(x_{b(i+1)})}{b},\ i=0,1,\cdots.$$
We denote it by $H_{\triangle,u}^{\infty}$.
\end{definition}

The twisted Hodge  polygon for Laurent polynomials can be found in
the literature, see Adolphson-Sperber [AS,AS2]. Liu \cite{Liu2}
proved the following.

\begin{theorem} We have
$$T-adic \text{ NP of }
C_{f,\chi}(s,T)\geq  {\rm ord}_p(q)(p-1)H_{\triangle,u}^{\infty},$$
 where NP is the short for Newton polygon.
\end{theorem}

If $\triangle$ is dimension one, the Newton polygon of the
$L$-function $L_{f,\chi}(s,\pi_m)$  for $m=1$ was studied  by
Blache-F\'{e}rard-Zhu \cite{BFZ}, and if $\chi$ is trivial, it was
studied by Liu-Liu-Niu \cite{LLN} for $m\geq 1$.

From now on, we assume that $\triangle=[0,d]$, $\chi=\omega^{-u}$
with $1\leq u\leq q-1$, and $q=p^b$.  Write
$$u=u_0+u_1p+\cdots+u_{b-1}p^{b-1}, 0\leq u_i\leq p-1.$$

\begin{definition}For $a\in \mathbb{N}$, $1\leq i\leq b$,
$$\delta_{\in}^{(i)}(n)=\left\{
                        \begin{array}{ll}
                          1, & \hbox{}pl\equiv n-u_{b-i}(d)~~\text{ for some }~~ l<d\{\frac{n}{d}\}; \\
                          0, & \hbox{}otherwise.
                        \end{array}
                      \right.
$$where $\{\cdot\}$ is the fractional part of a real number.
\end{definition}
\begin{definition}
A convex function on $\mathbb{R}_{\geq0}$ which is linear between
consecutive integers with initial value $0$ is called the twisted
arithmetic polygon of $\triangle =[0,d]$ if  its slopes between
consecutive integers are the numbers
$$\omega_{\triangle,u}(n)=\frac{1}{b}\sum\limits_{i=1}^{b }(\lceil\frac{(p-1)n+u_{b-i}}{d}\rceil-\delta_{\in}^{(i)}(n)), a\in \mathbb{N},$$
where $\lceil\cdot\rceil$ is the least integer equal or greater than
a real number.  We denote it by $p_{\triangle,u}$.
\end{definition}

Liu-Niu\cite{LN} proved the following.
\begin{theorem}
If $p>4d$,
$$ T-adic \text{ NP of }
C_{f,\chi}(s,T)\geq bp_{\triangle,u}.$$
\end{theorem}
By a result of Li\cite{Li}, $L_{f,\chi }(s,\pi_m)$ is a polynomial
with degree $p^{m-1}d$ if $p\nmid d$. Combined this result with the
above theorem, one can infer the following.
\begin{theorem}If $p>4d$,  then
$$\pi_m-adic \text{ NP of } L_{f,\chi}(s,\pi_m)\geq bp_{\triangle,u}\text{ on }[0,p^{m-1}d],$$
with equality holding for a generic $f$ of degree $d$.
\end{theorem}

We assume that the second highest exponent of $f$ is $ k$. So $k
\leq d-1$, and
$$f(x)=(a_dx^d,0,0,\cdots )+\sum\limits_{i=1}^k(a_ix^i,0,0,\cdots )\in W(\mathbb{F}_q[x])\text{ with }a_da_k\neq
0.$$

For  $a\in  \mathbb{N}$, define

\begin{align*}
\varpi_{d,[0,k],u}(a)=&\frac{1}{b}\sum\limits_{i=1}^{b}([\frac{pa+u_{b-i}}{d}]-[\frac{a}{d}]+[\frac{r_{a,i}}{k}]
-[\frac{r_a}{k}])\\
&+\frac{1}{b}\sum\limits_{i=1}^{b }\sum\limits_{j=0}^{r_a}(
1_{\{\frac{r_{j,i}}{k}\}>\{\frac{r_a}{k}\}}-1_{\{\frac{r_j}{k}\}>\{\frac{r_a}{k}\}})\\
&-\frac{1}{b}\sum\limits_{i=1}^{b }\sum\limits_{j=0}^{r_{a-1}}(
1_{\{\frac{r_{j,i}}{k}\}>\{\frac{r_{a-1}}{k}\}}-
1_{\{\frac{r_j}{k}\}>\{\frac{r_{a-1}}{k}\}}),
\end{align*}
where  $r_a=d\{\frac{a}{d}\}, r_{a,i}=d\{\frac{pa+u_{b-i}}{d}\}$ for
$a\in \mathbb{N}, 1\leq i\leq b $.

\begin{definition}A convex function on
$\mathbb{R}_{\geq 0}$ which is linear between consecutive integers
with initial value $0$ is called the twisted arithmetic polygon of
$\{d\}\cup[0,k]$ if its slopes between consecutive integers are the
numbers $\varpi_{d,[0,k],u}(a)$, $a\in \mathbb{N}$. We denote it by
$p_{d,[0,k],u}$ .
\end{definition}
We can prove the following.
\begin{theorem} \label{upper-bound} We have $p_{d,[0,k],u}\geq p_{\triangle,u}$.
\end{theorem}
 The  main result of this paper is
the following.
\begin{theorem}\label{first}If $p>d(2d+1)$,  then $$T-adic \text{ NP of } C_{f}(s,T)\geq
 \text{ord}_p(q)p_{d,[0,k],u}.$$
\end{theorem}
\begin{corollary}\label{second} If $p>d(2d+1)$, then
  $$\pi_m\text{ -adic NP of }
L_{f,\chi}(s,\pi_m) \geq \text{ord}_p(q)p_{d,[0,k],u}  \text{ on }
[0,p^{m-1}d].$$
\end{corollary}

\section{The $T$-adic Dwork Theory}In this section we review the $T$-adic analogue of Dwork
theory on exponential sums.

We can write
$$\frac{u}{q-1}=-(u_0+u_1p+\cdots),~~~u_i=u_{b+i} \text{ for } i\geq 0.$$
and  $p^iu=q_i(q-1)+s_i$ for $i\in \mathbb{N}$ with $0\leq s_i<q-1,$
then $s_{b-l}=u_l+u_{l+1}p+\cdots+u_{b+l-1}p^{b-1}$ for $0\leq l\leq
b-1$ and $s_i=s_{b+i}.$

Write $C_u=\{v\in\mathbb{N}|v\equiv u({ \rm mod } q-1)\}$. Let
$$B_u=\{\sum_{v\in
C_u}b_v\pi^{\frac{v}{d(q-1)}}x^{\frac{v}{q-1}}:b_v\in\mathbb{Z}_q[[\pi^{\frac{1}{d(q-1)}}]]
\text{ and } {\rm ord}_{\pi}b_v\rightarrow\infty \text{ as }
v\rightarrow\infty\}.$$ Define
$B=\bigoplus\limits_{i=1}^{b}B_{p^iu},$ then
$B=\bigoplus\limits_{i=1}^{b}B_{p^iu}$ has a basis represented by
$$\coprod\limits_{1\leq i\leq
b}\{x^{\frac{s_i}{q-1}+j}\}_{j\in\mathbb{N}}.$$

 Note
that the Galois group of $\mathbb{Q}_q$ over $\mathbb{Q}_p$ can act
on $B$ but keeping $\pi^{1/d}$ as well as the variable $x$ fixed.
 Let $\sigma\in{\rm Gal}(\mathbb{Q}_q/\mathbb{Q}_p)$ be the Frobenius element  such that
$\sigma(\zeta)=\zeta^p$ if $\zeta$ is a $(q-1)$-th root of unity and
$\Psi_p$  the operator on $B$ defined by the formula
$$\Psi_p(
 \sum\limits_{i\in
\mathbb{N}} c_ix^i)=\sum\limits_{i\in \mathbb{N}} c_{pi}x^i.$$

  The Frobenius operator $\Psi$ on B is defined by
$\Psi:=\sigma^{-1}\circ\Psi_p\circ E_f$,
 where $E_f(x) :=E(\pi \hat{a}_dx^d)\prod\limits_{i=1}^{k}E(\pi
 \hat{a}_ix^i)$.

Note that $\Psi:B_u\rightarrow B_{p^{-1}u}=B_{p^{b-1}u},$ hence
$\Psi$ is well defined. It follows that $\Psi^b$ operates on $B_u$
and is linear over $\mathbb{Z}_q[[\pi^{\frac{1}{d(q-1)}}]].$
Moreover it is completely continuous in the sense of Serre
\cite{Se}.

\begin{theorem}[$T$-adic Dwork trace formula]\label{analytic-trace-formula}
$$C_{f,\chi}(s,T)=\det(1-\Psi^bs|B_u/\mathbb{Z}_q[[\pi^{\frac{1}{d(q-1)}}]]).$$
\end{theorem}

\section{Key estimate}
In order to study
$$C_{f,\chi}(s,T)=\det(1-\Psi^bs|B_u/\mathbb{Z}_q[[\pi^{\frac{1}{d(q-1)}}]]), $$
we first study
$$\det(1-\Psi
s|B/\mathbb{Z}_p[[\pi^{\frac{1}{d(q-1)}}]])=\sum_{i=0}^{\infty}(-1)^ic_is^i.$$
We are going to show that
\begin{theorem}  \label{finalreduction}If $p>d(2d+1),$ then we
have $${\rm ord}_{\pi}(c_{b^2m})\geq b^2p_{d,[0,k],u}(m).$$
\end{theorem}
Consider the operator $\Psi_p\circ E_f(x)$ on $B$,  we have
\begin{align*} \Psi_p\circ
E_f(x)(x^{\frac{s_i}{q-1}+j})&=\Psi_p(\sum_{l=0}^{\infty}\gamma_lx^{\frac{s_i}{q-1}+j+l})\\
&=\sum_{pl+s_i-j\geq0}\gamma_{pl+s_i-j}x^{l+\frac{p^{b-1}s_i}{q-1}}\\
&=\sum_{l=0}^{\infty}\gamma_{pl+u_{b-i}-j}x^{l+\frac{s_{i-1}}{q-1}}.
\end{align*}

Then the matrix of $\Psi_p\circ E_f(x)$ on B with respect to the
basis $\coprod\limits_{1\leq i\leq b
}\{x^{\frac{s_i}{q-1}+j}\}_{j\in\mathbb{N}}$ is
$$(G_{(k,l)(i,j)})_{ 1\leq k ,i\leq b , l,j\in \mathbb{N}}.$$
where $$G_{(k,l)(i,j)}=\left\{
                        \begin{array}{ll}
                        \gamma_{pl+u_{b-i}-j}, & \hbox{} k=i-1;\\
                          0, & \hbox{} otherwise.
                        \end{array}
                      \right.
$$

Fix a normal basis $\bar{\xi}_1,\cdots,\bar{\xi}_{b }$ of
$\mathbb{F}_q$ over $\mathbb{F}_p$. Let $\xi_1,\cdots,\xi_{b}$ be
their Teichm\"{u}ller lifts. Then $\xi_1,\cdots,\xi_{b }$ is a
normal basis of $\mathbb{Q}_q$ over $\mathbb{Q}_p$, and $\sigma$
acts on $\xi_1,\cdots,\xi_{b }$ as a permutation.  Write
$$\xi_v^{\sigma^{-1}}G_{(k,l)(i,j)}^{\sigma^{-1}}=\sum_{w=1}^{b }G_{((k,l),w)((i,j),v)}\xi_w.$$

It is easy to see that $G_{((k,l),w)((i,j),v)}=0$ if $k\neq i-1.$
For $k=i-1,$ write
$$G_{((i-1,l),w)((i,j),v)}=G_{(l,w)(j,v)}^{(i)}.$$
\\

Write $G^{(i)}=(G_{(l,w)(j,v)}^{(i)})_{l,j\in \mathbb{N}, 1\leq w
,v\leq b },$ then the matrix of the operator $\Psi$ on B over
$\mathbb{Z}_p[[\pi^{\frac{1}{d(q-1)}}]]$ with respect to the basis
$\{\xi_vx^{\frac{s_i}{q-1}+j}\}_{1\leq i,v\leq b;j\in\mathbb{N}}$ is

$$G=\left(
                \begin{array}{ccccc}
                  0 & G^{(1)} & 0 & \cdots & 0 \\
                  0 & 0 & G^{(2)} & \cdots & 0 \\
                  \vdots & \vdots & \vdots & \vdots & \vdots \\
                  0 & 0 & 0 & \cdots & G^{(b-1)} \\
                  G^{(b)}& 0 & 0 & \cdots & 0 \\
                \end{array}
              \right).
$$

Hence we have $$\det(1-\Psi
s|B/\mathbb{Z}_p[[\pi^{\frac{1}{d(q-1)}}]])=\det(1-G
s)=\sum_{m=0}^{\infty}(-1)^mc_ms^m,$$ with
$c_m=\sum\limits_{F}\det(F)$, where $F$ runs over all principle
$m\times m$ submatrix  of $G.$

For every principle submatrix $F$ of $G,$ write $F^{(i)}=F\cap
G^{(i)}$ as the submatrix of $G^{(i)}.$ For principle $bm\times bm$
submatrix $F$ of $G,$ by linear algebra,   if one of $F^{(i)}$ is
not $m\times m$ submatrix of $G^{(i)}$ , then at least one row or
column of $F$ are 0 since $F$ is principle.

Let $\mathcal{F}_m$ be the set of all $bm\times bm$ principle
submatrices F of $G$ with $F^{(i)}$ all $m\times m$ submatrices of
$G^{(i)}$ for each $1\leq i\leq b.$

\begin{lemma}We have
$$c_{bm}=\sum_{F\in\mathcal{F}_m}\det(F)=\sum_{F\in\mathcal{F}_m}(-1)^{m^2(b-1)}\prod_{i=1}^b\det(F^{(i)}).$$
\end{lemma}

\begin{corollary}\label{Cbn}
$$c_{b^2m}=\sum_{F\in\mathcal{F}_{bm}}\prod_{i=1}^b\det(F^{(i)})=\sum_{F\in\mathcal{F}_{bm}}
\prod_{i=1}^b(\sum\limits_{\tau}sgn(\tau)\prod\limits_{
(l,\omega)\in R_i}G_{(l,w)\tau(l,w)}^{(i)}),$$

where $R_i$  runs over all subsets of $
\mathbb{N}\times\{1,2,\cdots,b\}$ with cardinality $bm$, $\tau$ runs
over all permutations of $R_i$, $1\leq i\leq b$.
\end{corollary}

So Theorem \ref{finalreduction} is reduced to the following.
\begin{theorem}\label{firsthalf} Let $p>d(2d+1)$. Then we have
$$\sum_{i=1}^b{\rm ord}_{\pi}(\sum\limits_{\tau}sgn(\tau)\prod\limits_{
(l,\omega)\in R_i}G_{(l,w)\tau(l,w)}^{(i)})\geq
b^2p_{d,[0,k],u}(m).$$
\end{theorem}

Let $O(\pi^{\alpha})$ denote any element of $\pi$-adic order $\geq
\alpha$.
\begin{lemma}
For any $1\leq i\leq b$ and $1\leq w,v\leq b,$ we have
$$G_{(l,w)(j,v)}^{(i)}=O(\pi^{[\frac{pl+u_{b-i}-j}{d}]+\lceil \frac{d\{\frac{pl+u_{b-i}-j}{d}\}}{k}\rceil}).$$
\end{lemma}

\proof Write $$E_f(x)= \sum\limits_{i\in \mathbb{N}}\gamma_ix^i.$$
 Then $$\gamma_i=\sum\limits_{\stackrel{dn_d+\sum\limits_{j=1}^kjn_j=i}
{n_j\geq0}}\pi^{\sum\limits_{j=1}^kn_j+n_d}
\prod\limits_{j=1}^k\lambda_{n_j}\hat{a}_j^{n_j}\lambda_{n_d}\hat{a}_d^{n_d}=O(\pi^{[\frac{i}{d}]+\lceil\frac{r_i}{k}\rceil}).$$
Since
$$\xi_v^{\sigma^{-1}}G_{(i-1,l)(i,j)}^{\sigma^{-1}}=\xi_v^{\sigma^{-1}}\gamma_{pl+u_{b-i}-j}^{\sigma^{-1}}
=\sum_{w=1}^{b}G_{(l,w)(j,v)}^{(i)}\xi_w, $$we have $${\rm
ord}_{\pi}(G_{(l,w)(j,v)}^{(i)})={\rm
ord}_{\pi}(\gamma_{pl+u_{b-i}-j}).$$The lemma now follows.\endproof

By the above lemma, Theorem \ref{firsthalf} is reduced to the
following.
\begin{theorem}
Let $p>d(2d+1)$.  For  $1\leq i \leq b,$ let $R_i\subset
\mathbb{N}\times\{1,2,\cdots,b\}$ be a subset of cardinality $bm$,
and $\tau$ a permutation of $R_i$.  Then
$$\sum_{i=1}^b\sum\limits_{(l,\omega)\in
R_i}([\frac{pl+u_{b-i}-\tau(l)}{d}]+\lceil\frac{d\{\frac{pl+u_{b-i}-\tau(l)}{d}\}}{k}\rceil)
\geq b^2p_{d,[0,k],u}(m),$$ where $\tau(l)$ is defined by
$\tau(l,\omega)=(\tau(l), \tau(\omega))$.
\end{theorem}

\proof  By definition, we have
$$p_{d,[0,k],u}(m)=\sum_{a=0}^{m-1}\varpi_{d,[0,k],u}(a)=
\sum_{i=1}^bp_{d,[0,k],u}^{(i)}(m),
$$
where
$$p_{d,[0,k],u}^{(i)}(m)=\frac{1}{b}\sum_{a=0}^{m-1}\bigg([\frac{pa+u_{b-i}}{d}]-[\frac{a}{d}]+[\frac{r_{a,i}}{k}]-[\frac{r_a}{k}]+
1_{\{\frac{r_{a,i}}{k}\}>\{\frac{r_{m-1}}{k}\}}-
1_{\{\frac{r_a}{k}\}>\{\frac{r_{m-1}}{k}\}}\bigg).
$$

Then it suffices to show that for $1\le i \le b$,  and for any
permutation $\tau$ of $R_i$, we have

$$\sum\limits_{(l,\omega)\in
R_i}([\frac{pl+u_{b-i}-\tau(l)}{d}]+\lceil\frac{d\{\frac{pl+u_{b-i}-\tau(l)}{d}\}}{k}\rceil)\geq
b^2p_{d,[0,k],u}^{(i)}(m).$$

Note that $$\sum\limits_{(l,\omega)\in
R_i}([\frac{pl+u_{b-i}-\tau(l)}{d}]+\lceil\frac{d\{\frac{pl+u_{b-i}-\tau(l)}{d}\}}{k}\rceil)$$
$$=\sum\limits_{(l,\omega)\in R_i}
([\frac{pl+u_{b-i}}{d}]-[\frac{l}{d}]+[\{\frac{pl+u_{b-i}}{d}\}-\{\frac{\tau(l)}{d}\}]+\lceil\frac{r_{l,i}-r_{\tau(l)}
-d[\{\frac{pi}{d}\}-\{\frac{\tau(i)}{d}\}]}{k}\rceil)$$
$$=\sum\limits_{(l,\omega)\in
R_i}([\frac{pl+u_{b-i}}{d}]-[\frac{l}{d}]+[\frac{r_{l,i}}{k}]-[\frac{r_l}{k}]+\lceil\frac{(d-k)1_{r_{\tau(l)}>r_{l,i}}}{k}
+\{\frac{r_{l,i}}{k}\}-\{\frac{r_{\tau(l)}}{k}\}\rceil).$$
 And
$$\sum\limits_{(l,\omega)\in
R_i}(\lceil\frac{(d-k)1_{r_{\tau(l)}>r_{l,i}}}{k}
+\{\frac{r_{l,i}}{k}\}-\{\frac{r_{\tau(l)}}{k}\}\rceil)$$
$$\geq
\sum\limits_{(l,\omega)\in R_i}1_{\{\frac{r_{\tau(l) }}{k}\}\leq
\{\frac{r_{m-1} }{k}\}<\{\frac{r_{l,i}}{k}\}}$$ $$\geq
\sum\limits_{(l,\omega)\in
R_i}1_{\{\frac{r_{l,i}}{k}\}>\{\frac{r_{m-1}}{k}\}}-
\sum\limits_{(l,\omega)\in R_i}
1_{\{\frac{r_l}{k}\}>\{\frac{r_{m-1}}{k}\}}.
$$

We have

$$\sum\limits_{(l,\omega)\in
R_i}([\frac{pl+u_{b-i}}{d}]-[\frac{l}{d}]+[\frac{r_{l,i}}{k}]-[\frac{r_l}{k}])$$
$$=b\sum\limits_{l=0}^{m-1}([\frac{pl+u_{b-i}}{d}]-[\frac{l}{d}]+[\frac{r_{l,i}}{k}]-[\frac{r_l}{k}])
+\sum\limits_{\stackrel{(l,\omega)\in R_i}{l\geq
m}}([\frac{pl+u_{b-i}}{d}]-[\frac{l}{d}]+[\frac{r_{l,i}}{k}]-[\frac{r_l}{k}])$$
$$\quad -\sum\limits_{\stackrel{(l,\omega)\notin R_i}{0\leq l<m}}
([\frac{pl+u_{b-i}}{d}]-[\frac{l}{d}]+[\frac{r_{l,i}}{k}]-[\frac{r_l}{k}])$$
$$\geq
b\sum\limits_{l=0}^{m-1}([\frac{pl+u_{b-i}}{d}]-[\frac{l}{d}]+[\frac{r_{l,i}}{k}]-[\frac{r_l}{k}])
+N([\frac{pm+u_{b-i}}{d}]-[\frac{m}{d}]-[\frac{d-1}{k}])$$
$$-N([\frac{p(m-1)+u_{b-i}}{d}]-[\frac{m-1}{d}]+[\frac{d-1}{k}])$$
$$\geq
b\sum\limits_{l=0}^{m-1}([\frac{pl+u_{b-i}}{d}]-[\frac{l}{d}]+[\frac{r_{l,i}}{k}]-[\frac{r_l}{k}])
+N([\frac{p}{d}]-1-2(d-1)),
$$
where $N=\#\{(l,\omega)\in R_i| l\geq m\}=\#\{(l,\omega)\notin
R_i|0\leq l <m\}$.

 Similarly, we have

$$\sum\limits_{(l,\omega)\in
R_i}1_{\{\frac{r_{l,i}}{k}\}>\{\frac{r_{m-1}}{k}\}} \geq
b\sum\limits_{l=0}^{m-1}1_{\{\frac{r_{l,i}}{k}\}>\{\frac{r_{m-1}}{k}\}}-N,
$$ and
$$
\sum\limits_{(l,\omega)\in R}
1_{\{\frac{r_l}{k}\}>\{\frac{r_{m-1}}{k}\}}\leq
b\sum\limits_{l=0}^{m-1}1_{\{\frac{r_l}{k}\}>\{\frac{r_{m-1}}{k}\}}+N.
$$
Therefore, $$\sum\limits_{(l,\omega)\in
R_i}([\frac{pl+u_{b-i}-\tau(l)}{d}]+\lceil\frac{d\{\frac{pl+u_{b-i}-\tau(l)}{d}\}}{k}\rceil)$$
$$\geq
b\sum\limits_{l=0}^{m-1}([\frac{pl+u_{b-i}}{d}]-[\frac{l}{d}]+[\frac{r_{l,i}}{k}]-[\frac{r_l}{k}]+1_{\{\frac{r_{l,i}}{k}\}>\{\frac{r_{m-1}}{k}\}}
-1_{\{\frac{r_l}{k}\}>\{\frac{r_{m-1}}{k}\}})$$
$$+N([\frac{p}{d}]-1-2(d-1)-2)$$$$=b^2p_{d,[0,k],u}^{(i)}(m)+N([\frac{p}{d}]-2d+1)\geq
b^2p_{d,[0,k],u}^{(i)}(m).
$$
\endproof
\section{Proof of the main result}
In this section we prove Theorem \ref{first} , which says that, if
$p>d(2d+1)$,  then $$T-adic \text{ NP of } C_{f,\chi}(s,T)\geq
bp_{d,[0,k],u}.$$

\begin{lemma}\label{zqtozp} We have
$$T-adic \text{ NP of }
 \det(1-\Psi^bs^b|B/\mathbb{Z}_q[[\pi^{\frac{1}{d(q-1)}}]])$$
 $$=T-adic \text{ NP of }
\det(1-\Psi s|B/\mathbb{Z}_p[[\pi^{\frac{1}{d(q-1)}}]]).$$
\end{lemma}

\proof The lemma follows from the following:
\begin{align*}
\prod_{\zeta^b=1}\det(1-\Psi\zeta
s|B/\mathbb{Z}_p[[\pi^{\frac{1}{d(q-1)}}]])&=\det(1-\Psi^bs^b|B/\mathbb{Z}_p[[\pi^{\frac{1}{d(q-1)}}]])\\
&={\rm
Norm}(\det(1-\Psi^bs^b|B/\mathbb{Z}_q[[\pi^{\frac{1}{d(q-1)}}]])),
\end{align*}
where the Norm is the norm map from
$\mathbb{Q}_q[[\pi^{\frac{1}{d(q-1)}}]]$ to
$\mathbb{Q}_p[[\pi^{\frac{1}{d(q-1)}}]].$ \qed

\begin{lemma}\label{BtoB_i} We have
$$T-adic \text{ NP of }C_{f,\chi}(s,T)^b =T-adic \text{ NP of }
of \det(1-\Psi^bs|B/\mathbb{Z}_q[[\pi^{\frac{1}{d(q-1)}}]]).$$
\end{lemma}

\proof Let $\sigma$ act on $\mathbb{Q}_q[[T]]$ coordinate-wise.
Hence \begin{align*}
S_{f,\chi^p}(l,T)&=\sum_{x\in\mathbb{F}_{q^l}^{\times}}\chi({\rm
Norm}_{\mathbb{F}_{q^l}/\mathbb{F}_q}(x))^p(1+T)^{Tr_{\mathbb{Q}_{q^l}/\mathbb{Q}_p}(f(x))}\\
&=S_{f,\chi}(l,T)^{\sigma},
\end{align*}
therefore $C_{f,\chi^p}(s,T)=C_{f,\chi}(s,T)^{\sigma},$ which yields
that the T-adic Newton polygons of $C_{f,\chi^p}(s,T)$ and
$C_{f,\chi}(s,T)$ coincide with each other. Hence the lemma follows
from the following
\begin{align*}
\prod_{i=1}^{b}C_{f,\chi}(s,T)^{\sigma^i}&=\prod_{i=1}^{b}C_{f,\chi^{p^i}}(s,T)\\
&=\prod_{i=1}^{b}\det(1-\Psi^bs|B_{p^iu}/\mathbb{Z}_q[[\pi^{\frac{1}{d(q-1)}}]])\\
&=\det(1-\Psi^bs|B/\mathbb{Z}_q[[\pi^{\frac{1}{d(q-1)}}]]).
\end{align*}
\qed

\begin{corollary}\label{Npoints}
The $T$-adic Newton polygon of $C_{f,\chi}(s,T)$ is the lower convex
closure of the points $$(i,\frac{1}{b}{\rm
ord}_{\pi}c_{b^2i}),~~~i=0,1,\cdots$$
\end{corollary}
\proof By Lemma \ref{zqtozp}, the T-adic Newton polygon of
$$\det(1-\Psi^bs^b|B/\mathbb{Z}_q[[\pi^{\frac{1}{d(q-1)}}]])$$ is
the lower convex closure of the points $$(bi,{\rm
ord}_{\pi}c_{bi}),~~~i=0,1,\cdots.$$ Hence the T-adic Newton polygon
of $$\det(1-\Psi^bs|B/\mathbb{Z}_q[[\pi^{\frac{1}{d(q-1)}}]])$$ is
the convex closure of the points $$(i,{\rm
ord}_{\pi}c_{bi}),~~~i=0,1,\cdots.$$ By Lemma \ref{BtoB_i}, the
T-adic Newton polygon of $C_{f,\chi}(s,T)^b$ is the lower convex
closure of the points $$(bi,{\rm
ord}_{\pi}c_{b^2i}),~~~i=0,1,\cdots.$$ The lemma is proved. \qed
\\

We now prove  Theorem \ref{first}.

{\it Proof of Theorem \ref{first}. }  By Theorem
\ref{analytic-trace-formula}, we have
$$C_{f,\chi}(s,T)=\det(1-\Psi^bs|B_u/\mathbb{Z}_q[[\pi^{\frac{1}{d(q-1)}}]]),
$$
Then by Corollary \ref{Npoints}, the $T$-adic Newton polygon of
$C_{f,\chi}(s,T)$ is the lower convex closure of the points
$$(i,\frac{1}{b}{\rm ord}_{\pi}c_{b^2i}),~~~i=0,1,\cdots.$$ Therefore
the result follows from Theorem \ref{finalreduction}, which says
that, if $p>d(2d+1),$ then we have
$${\rm ord}_{\pi}(c_{b^2m})\geq b^2p_{d,[0,k],u}(m).$$\qed

We conclude this section by proving Corollary \ref{second}.

 {\it
Proof of Corollary \ref{second}. } Assume that $
L_{f,\chi}(s,\pi_m)=\prod\limits_{i=1}^{p^{m-1}d}(1-\beta_is) $.
Then
$$C_{f,\chi}(s,\pi_m )= \prod\limits_{j=0}^{\infty} L_{f,\chi}(q^js, \pi_m)=
\prod\limits_{j=0}^{\infty}\prod\limits_{i=1}^{p^{m-1}d}(1-\beta_iq^js).$$
Therefore the slopes of the  $q$-adic Newton polygon of
$C_{f,\chi}(s,\pi_m)$ are the numbers
$$j+{\rm ord}_q(\beta_i),\
1\leq i\leq p^{m-1}d, j=0,1,\cdots.$$ It is well-known that ${\rm
ord}_q(\beta_i)\leq 1$ for all $i$. Therefore, $$q-adic \text{ NP of
}L_{f,\chi}(s,\pi_m)= q-adic\text{ NP of }C_{f,\chi}(s,\pi_m) \text
{ on } [0,p^{m-1}d].$$ It follows that
 $$\pi_m-adic \text{ NP of
}L_{f,\chi}(s,\pi_m)= \pi_m-adic\text{ NP of }C_{f,\chi}(s,\pi_m)
\text { on } [0,p^{m-1}d].$$ By  the integrality of
$C_{f,\chi}(s,T)$ and Theorem \ref{first}, we have
$$\pi_m-adic \text{ NP of } C_{f,\chi}(s,\pi_m)\geq T-adic \text{ NP of } C_{f,\chi}(s,T)\geq \text{ord}_p(q)p_{d,[0,k],u}.$$

Therefore,   $$\pi_m\text{ -adic NP of } L_{f,\chi}(s,\pi_m) \geq
\text{ord}_p(q)p_{d,[0,k],u}  \text{ on } [0,p^{m-1}d].$$ \qed

\section{comparison between arithmetic polygons}
In this section we prove Theorem \ref{upper-bound}, which says that
$$p_{d,[0,k],u}\geq p_{\triangle,u}.$$

{\it Proof of Theorem \ref{upper-bound} }  It is clear that
$p_{d,[0,k]}(0)=p_{\triangle}(0)$.   It suffices to show that, for
$m\in \mathbb{N}$, we have $$p_{d,[0,k],u}(m+1) \geq
p_{\triangle,u}(m+1).$$ By a result in Liu-Niu \cite{LN}, we have
$$p_{\triangle,u}(m+1)=\sum_{i=1}^bp_{\triangle,u}^{(i)}(m+1),
$$
where
$$p_{\triangle,u}^{(i)}(m+1)=\frac{1}{b}\bigg(\sum_{a=0}^{m }(\lceil\frac{pa+u_{b-i}}{d}\rceil-
\lceil\frac{a}{d}\rceil)+\sum\limits_{a=0}^{r_m}(1_{\{\frac{a}{d}\}^{'}\leq
\frac{r_m}{d}<\{\frac{pa+u_{b-i}}{d}\}^{'}}-1_{\{\frac{u_{b-i}}{d}\}^{'}\leq
\frac{r_m}{d}})\bigg).$$

By definition, we have
$$p_{d,[0,k],u}(m+1)=\sum_{a=0}^{m }\varpi_{d,[0,k],u}(a)=
\sum_{i=1}^bp_{d,[0,k],u}^{(i)}(m+1),
$$
where
$$p_{d,[0,k],u}^{(i)}(m+1)=\frac{1}{b}\sum_{a=0}^{m }\bigg([\frac{pa+u_{b-i}}{d}]-[\frac{a}{d}]+[\frac{r_{a,i}}{k}]-[\frac{r_a}{k}]+
1_{\{\frac{r_{a,i}}{k}\}>\{\frac{r_{m }}{k}\}}-
1_{\{\frac{r_a}{k}\}>\{\frac{r_{m }}{k}\}}\bigg).
$$
Then it suffices to show that
$$p_{d,[0,k],u}^{(i)}(m+1)\geq p_{\triangle,u}^{(i)}(m+1).$$

 \noindent For $m\geq 0$, $1\leq i \leq b$,
let $A_{i1}=\{0\leq a\leq r_m|a\neq r_{l,i} \text{ for some } 0\leq
l\leq r_m\}$. Note that
$$\{a:0\leq a\leq r_m\} =A_{i1}\cup A_{i2},$$
 where $A_{i2}= \{r_{ a,i}|0 \leq a, r_{a,i}\leq r_m\}.$
And
$$\{r_{ a,i }|0\leq a\leq r_m\}=A_{i2}\cup A_{i3} ,$$
where $A_{i3}=\{r_{ a,i}>r_m| 0 \leq a\leq r_m\}$, so we have
$|A_{i1}|=|A_{i3}|$.   Then

$$\sum\limits_{a=0}^{m }([\frac{r_{a,i}}{k}]-[\frac{r_a}{k}]+
1_{\{\frac{r_{a,i}}{k}\}>\{\frac{r_m}{k}\}}-
1_{\{\frac{r_a}{k}\}>\{\frac{r_m}{k}\}})$$
$$=\sum\limits_{a=0}^{r_m}([\frac{r_{a,i}}{k}]-[\frac{a}{k}]+
 1_{\{\frac{r_{a,i}}{k}\}>\{\frac{r_m}{k}\}}-
 1_{\{\frac{a}{k}\}>\{\frac{r_m}{k}\}})$$
$$=\sum\limits_{r_{a,i}\in A_{i3}}[\frac{r_{a,i}}{k}]-\sum\limits_{a\in
A_{i1}}[\frac{a}{k}]+ \sum\limits_{r_{a,i}\in
A_{i3}}1_{\{\frac{r_{a,i}}{k}\}>\{\frac{r_m}{k}\}}-
\sum\limits_{a\in A_{i1}} 1_{\{\frac{a}{k}\}>\{\frac{r_m }{k}\}}$$
$$=\sum\limits_{r_{a,i}\in A_{i3}}\lceil\frac{r_{a,i}-r_m}{k}\rceil+
\sum\limits_{a\in
A_{i1}}(\lceil\frac{r_m-a}{k}\rceil-1_{\{\frac{a}{k}\}>\{\frac{r_m}{k}\}}-1_{\{\frac{a}{k}\}<\{\frac{r_m}{k}\}})$$
$$\geq \sum\limits_{r_{a,i}\in A_{i3}}\lceil\frac{r_{a,i}-r_m}{k}\rceil.$$

Therefore,   we have
$$b(p_{d,[0,k],u}^{(i)}(m+1)-p_{\triangle,u}^{(i)}(m+1))$$
$$=\sum_{a=0}^{m
}\bigg([\frac{pa+u_{b-i}}{d}]-[\frac{a}{d}]+[\frac{r_{a,i}}{k}]-[\frac{r_a}{k}]+
1_{\{\frac{r_{a,i}}{k}\}>\{\frac{r_{m }}{k}\}}-
1_{\{\frac{r_a}{k}\}>\{\frac{r_{m }}{k}\}}\bigg)$$
$$-\bigg(\sum_{a=0}^{m }(\lceil\frac{pa+u_{b-i}}{d}\rceil-
\lceil\frac{a}{d}\rceil)+\sum\limits_{a=0}^{r_m}(1_{\{\frac{a}{d}\}^{'}\leq
\frac{r_m}{d}<\{\frac{pa+u_{b-i}}{d}\}^{'}}-1_{\{\frac{u_{b-i}}{d}\}^{'}\leq
\frac{r_m}{d}})\bigg)$$
$$\geq\sum_{a=0}^{m}(\lceil\{\frac{a}{d}\}\rceil-\lceil\{\frac{r_{a,i}}{d}\}\rceil)+\sum\limits_{r_{a,i}\in
A_{i3}} \lceil\frac{r_{a,i}-r_m
}{k}\rceil-\sum\limits_{a=1}^{r_m}1_{\frac{r_m}{d}<\{\frac{pa+u_{b-i}}{d}\}^{'}}+1_{\{\frac{u_{b-i}}{d}\}^{'}\leq
\frac{r_m}{d}}$$
$$=\sum_{a=0}^{r_m}(\lceil \frac{a}{d} \rceil-\lceil\{\frac{r_{a,i}}{d}\}\rceil)+\sum\limits_{a=1}^{r_m}
(\lceil\frac{r_{a,i}-r_m}{k}\rceil-1)1_{r_{a,i}>r_m}+\lceil\frac{r_{0,i}-r_m}{k}\rceil1_{\{\frac{u_{b-i}}{d}\}>\frac{r_m}{d}}$$
$$+\delta_{m,2}^{(i)}+1_{\{\frac{u_{b-i}}{d}\}^{'}\leq
\frac{r_m}{d}}$$
$$\geq
\delta_{m,1}^{(i)}+\delta_{m,2}^{(i)}+1_{\{\frac{u_{b-i}}{d}\}>\frac{r_m}{d}}+1_{\{\frac{u_{b-i}}{d}\}^{'}\leq
\frac{r_m}{d}}\geq 0.$$

where for $1\leq i\leq b$,
 $$ \delta_{m,1}^{(i)}=\left\{
                                  \begin{array}{ll}
                                    0, & \hbox{ if there exists }  0\leq a\leq r_m\text{ such that } pa+u_{b-i}\equiv 0(\text{ mod }d) ; \\
                                    -1, & \hbox{ otherwise }.
                                  \end{array}
                                \right.
$$
$$ \delta_{m,2}^{(i)}=\left\{
                                  \begin{array}{ll}
                                    -1, & \hbox{} \text{ if there exists } 1\leq a\leq r_m\text{ such that } pa+u_{b-i}\equiv 0(\text{ mod }d) ; \\
                                    0, & \hbox{ otherwise }.
                                  \end{array}
                                \right.
$$

The  theorem now follows. \qed

\end{document}